\documentclass[10pt,a4paper]{amsart}

\usepackage[a4paper, margin=1in]{geometry}
\usepackage[english,frenchb]{babel}
\usepackage[colorlinks=true]{hyperref}
\hypersetup{citecolor=blue}
\theoremstyle{plain}
\newtheorem{thrm}{\thrmname}[section]
\newtheorem{lmm}[thrm]{\lmname}

\newtheorem{prpstn}[thrm]{\prpstnname}
\providecommand{\thrmname}{Theorem}
\providecommand{\lmname}{Lemma}
\providecommand{\crllrname}{Corollary}
\providecommand{\prpstnname}{Proposition}
\theoremstyle{definition}
\newtheorem{dfntn}[thrm]{\dfntnname}

\newtheorem{rmrk}[thrm]{\rmrkname}
\providecommand{\dfntnname}{Definition}
\providecommand{\cnjctrname}{Conjecture}
\providecommand{\xmplname}{Example}
\providecommand{\rmrkname}{Remark}
\def\xC{\mathbb{C}} \def\xR{\mathbb{R}}
  \def\xN{\mathbb{N}}

\def\xCinfty{{\rm C}^{\infty}} 

\def\xLtwo{{\rm L}^{2}} 
 
\def\xLn#1{{\rm L}^#1}
\def\xdif{\,{\rm d}}
\newtheorem{ssmptn}[thrm]{Assumption}

\usepackage{xcolor}

\makeatletter
\newenvironment{abstracts}{%
  \ifx\maketitle\relax
    \ClassWarning{\@classname}{Abstract should precede
      \protect\maketitle\space in AMS document classes; reported}%
  \fi
  \global\setbox\abstractbox=\vtop \bgroup
    \normalfont\Small
    \list{}{\labelwidth\z@
      \leftmargin3pc \rightmargin\leftmargin
      \listparindent\normalparindent \itemindent\z@
      \parsep\z@ \@plus\p@
      
      \itemsep\medskipamount
    }%
}{%
  \endlist\egroup
  \ifx\@setabstract\relax \@setabstracta \fi
}

\newcommand{\abstractin}[1]{%
  \otherlanguage{#1}%
  \item[\hskip\labelsep\scshape\abstractname.]%
}
\makeatother

\begin{document}

%%-----------------------------
%%      the top matter
%%-----------------------------
\title[Smoothing effect for the Schrödinger equation]{Smoothing effect and quantum-classical correspondence\\for the Schrödinger equation with confining potential}%\thanks{...}\thanks{...}% At most 5 thanks
\author{Antoine Prouff}
\date{\today}
\address{Department of Mathematics, Purdue University, West Lafayette, IN, USA}
\email{aprouff@purdue.edu}
%\author{...}\address{...}
%\author{...}\address{...}
%
%\dedicated{\it Dedicated to Maurice Dupont} %if necessary
%

\begin{abstracts}
\abstractin{english}
The smoothing effect states that solutions to the Schrödinger equation in the Euclidean space have, for almost-every time, a local-in-space improved regularity (gain of half a derivative in Sobolev spaces). In this note, we show that, for the Schrödinger equation with a sub-quadratic confining potential, the smoothing effect is equivalent to an escape rate estimate on the associated classical flow. The proof relies on an Egorov theorem proved in~\cite{P:24Egorovinprep}.

\abstractin{french}
L'effet régularisant énonce que les solutions de l'équation de Schrödinger dans l'espace Euclidien ont, pour presque tout temps, un gain de régularité local en espace (gain d'une demie dérivée dans les espaces de Sobolev). Dans cette note, on montre, pour l'équation de Schrödinger avec potentiel confinant sous-quadratique, que l'effet régularisant est équivalent à une estimée du taux de fuite des trajectoires de la mécanique classique sous-jacente. La preuve repose sur un théorème d'Egorov tiré de~\cite{P:24Egorovinprep}.
\end{abstracts}
\selectlanguage{english}
\maketitle
%%-----------------------------
%%      your text
%%-----------------------------

\section{Introduction}

The simplest instance of smoothing effect for the Schrödinger equation, also know as Kato smoothing effect or {$1/2$-smoothing} effect, is the following: consider the Schrödinger equation in $\xR^d$, $d \ge 1$,
\begin{equation} \label{eq:SchLap}
i \partial_t \psi
	= \Delta \psi
\end{equation}
and denote by $t \mapsto e^{- i t \Delta} u$ the unique solution issued from $\psi(0, \bullet) = u \in \xLtwo(\xR^d)$. Then it holds
\begin{equation} \label{eq:1/2smoothing}
\forall \nu > 1/2, \exists C > 0 : \forall u \in \xLtwo(\xR^d) , \,
	 \int_\xR \left\| \left\langle x \right\rangle^{- \nu} (1 - \Delta)^{1/4} e^{- i t \Delta} u\right\|_{\xLtwo(\xR^d)}^2 \xdif t
		\le C \left\|u\right\|_{\xLtwo(\xR^d)}^2 ,
\end{equation}
where $\langle \bullet \rangle = \sqrt{1 + | \bullet |^2}$. In particular, one has $e^{- i t \Delta} u \in \xLtwo(\xR, {{\textrm{H}}_{\textrm{loc}}^{1/2}}(\xR^d))$.  This is a purely {$\xLtwo$-based} estimate, as opposed to Strichartz estimates, that quantify dispersion on the scale of Lebesgue spaces $\xLn{p}$, $p \in [1, \infty]$. We refer to lecture notes by Yajima~\cite{Yajima} and to the survey paper of Robbiano~\cite{Robbiano} for a thorough exposition of the Kato smoothing effect, its variants and generalizations, and an exhaustive list of references. Although the smoothing effect could be seen as an artifact of dispersion, variations of~\eqref{eq:SchLap} are available in a greater generality, including for operators having eigenvalues, in place of the Laplacian. One of the most striking results was given by Doi~\cite{Doi}, who proves an estimate similar to~\eqref{eq:1/2smoothing} with a Laplace-Beltrami operator associated to a non-trapping metric on $\xR^d$ as well as smooth magnetic and electric potentials with a sub-linear and sub-quadratic control at infinity respectively.

\subsection{Smoothing effect for Schrödinger operators with confining potential}

In this note, we consider the Schrödinger equation with a sub-quadratic potential:
\begin{equation} \label{eq:Sch}
i \partial_t \psi
	= P \psi ,
		\qquad
	P = - \tfrac{1}{2} \Delta + V(x) ,
\end{equation}
still in the {$d$-dimensional} Euclidean space.

\begin{ssmptn} \label{assum:V}
The potential $V$ is of class $\xCinfty$, non-negative, and satisfies for some $m \in (0, 1]$:
\begin{equation} \label{eq:assumV}
\begin{split}
\exists C > 0 : \forall x \in \xR^d , \qquad
	\dfrac{1}{C} \left\langle x \right\rangle^{2m} - C
		&\le V(x)
		\le C \left\langle x \right\rangle^{2m} \\
\forall \alpha \in \xN^d, \exists C_\alpha > 0 : \forall x \in \xR^d , \qquad
	\left|\partial^\alpha V(x)\right|
		&\le C_\alpha \left\langle x \right\rangle^{2 m - |\alpha|} .
\end{split}
\end{equation}
\end{ssmptn}

Denote by $t \mapsto e^{- i t P} u$ the unique solution to~\eqref{eq:Sch} issued from $\psi(0, \bullet) = u \in \xLtwo(\xR^d)$. In this setting, Doi's result reads as follows:

\begin{thrm}[Theorem 2.8 in {\cite{Doi}}] \label{thm:Doi}
For any $\nu > 1/2$ and all $T > 0$, there exists $\mathfrak{C}_0 > 0$ such that
\begin{equation} \label{eq:Doi}
\forall u \in \xLtwo(\xR^d) , \qquad
	\int_0^T \left\| \left\langle x \right\rangle^{-\nu} \mathrm{Op}\left( (1 + p)^{1/4} \right) e^{- i t P} u \right\|_{\xLtwo(\xR^d)}^2 \xdif t
		\le \mathfrak{C}_0 \left\|u\right\|_{\xLtwo(\xR^d)}^2 .
\end{equation}
Here $\mathrm{Op}((1 + p)^{1/4})$ is the Weyl quantization of the symbol $(1 + p)^{1/4}$ with $p(x, \xi) = \frac{1}{2} |\xi|^2 + V(x)$ (see Subsection~\ref{subsec:quantization}).
\end{thrm}

The method of Doi involves the construction of an escape function and a positive commutator argument, which were used later by Robbiano and Zuily~\cite{RZ} to generalize~\eqref{eq:Doi} to super-quadratic potentials (see also~\cite{RZexterior} for exterior domains). This strategy is generally understood as a way to formalize the fact that singularities propagate at infinite speed, and relies on the fact that the underlying metric is non-trapping (which is obviously the case for the flat Laplacian).

The purpose of this note is to recover Theorem~\ref{thm:Doi} as an explicit consequence of the quantum-classical correspondence principle. Our alternative proof relies on the so-called Egorov's theorem, which relates the quantum evolution ruled by~\eqref{eq:Sch} to the underlying classical dynamics naturally associated with our problem (see the following Subsection~\ref{subsec:classicaldynamics}). Although the method that we propose here seems to be less suitable for generalizations, it gives an interpretation of the smoothing effect as the quantum analogue to an estimate quantifying the ability of trajectories of the underlying classical dynamics to escape compact sets (see Proposition~\ref{prop:classicalsmoothing} and Remark~\ref{rem:classicalsmoothing} below).

Notice that under Assumption~\ref{assum:V}, the operator $P$ has compact resolvent, hence its spectrum is made of eigenvalues. Thus the word ``dispersion" is not really appropriate in this context; it is better to think of Theorem~\ref{thm:Doi} as a statement related in a quantitative way to the unbounded speed of propagation of energy. Beyond propagation of singularities, Egorov's theorem allows to describe propagation of energy. This more precise information allows to establish a parallel between the quantum smoothing estimate~\eqref{eq:Doi} and the classical dynamical estimate~\eqref{eq:classicalsmoothing} below, and clarifies the appearance of specific scales in~\eqref{eq:Doi} (namely the powers in $(1 + p)^{1/4}$ and $\langle x \rangle^{-1/2}$).

\subsection{The underlying classical dynamics} \label{subsec:classicaldynamics}

The operator $P = - \frac{1}{2} \Delta + V(x)$ can be seen as the quantization (see Subsection~\ref{subsec:quantization}) of the smooth function $p(x, \xi) = \frac{1}{2} |\xi|^2 + V(x)$ on the phase space $\xR^{2d} = \xR_x^d \times \xR_\xi^d$. The latter function is called the \emph{symbol} of $P$. The Hamiltonian flow $(\phi^t)_{t \in \xR}$ acting on $\xR^{2d}$ is the group of diffeomorphisms
\begin{equation*}
\phi^t : \xR^{2d} \longrightarrow \xR^{2d} ,
	\qquad
\phi^t \circ \phi^s = \phi^{t + s} , \;\, \forall t, s \in \xR ,
\end{equation*}
defined as $\phi^t(x, \xi) = (x^t, \xi^t)$, where
\begin{equation} \label{eq:Hamilton}
\dfrac{\xdif}{\xdif t} \begin{pmatrix}
x^t \\ \xi^t
\end{pmatrix}
	= \begin{pmatrix}
	\xi^t \\ - \nabla V\left(x^t\right)
	\end{pmatrix} .
\end{equation}
Thus $t \mapsto \phi^t(x, \xi)$ is the trajectory in phase space of a point mass subject to a force field $\vec F(x) = - \nabla V(x)$, given by Newton's second law of classical mechanics. This flow preserves the energy, namely
\begin{equation} \label{eq:conservationofenergy}
p \circ \phi^t = p , \qquad \forall t \in \xR .
\end{equation}

\subsection{Smoothing and escape rate of classical trajectories}

In~\cite{P:23}, we showed with elementary computations the following estimate at the level of classical mechanics.

\begin{prpstn}[Subsection 1.6 (1.34) in {\cite{P:23}}] \label{prop:classicalsmoothing}
Suppose $V$ is subject to Assumption~\ref{assum:V} and let $p(x, \xi) = \frac{1}{2} |\xi|^2 + V(x)$. Fix $T > 0$ and $\nu > 1/2$. Then there exists $\mathsf{C}_0 > 0$ such that
\begin{equation} \label{eq:classicalsmoothing}
\forall (x, \xi) \in \xR^{2d} , \qquad
	\sqrt{1 + p(x, \xi)} \int_0^T \dfrac{\xdif t}{\langle (\pi \circ \phi^t) (x, \xi) \rangle^{2 \nu}}
		\le \mathsf{C}_0 ,
\end{equation}
where $\pi : \xR^{2d} \to \xR^d$ is the projection onto the position variable $(x, \xi) \mapsto x$.
\end{prpstn}

This estimate was proved by exploiting the particular structure of the symbol $p$, and in particular the fact that $\Delta$ is the \emph{flat} Laplacian on $\xR^d$. It quantifies the ability of flow trajectories to escape compact sets. We claim that this is the classical counterpart to the smoothing inequality~\eqref{eq:Doi}. Such an inequality appears already in previous works as a quantitative version of the non-trapping condition. See for instance the paper of Burq~\cite[(1.4)]{Burq} on the Schrödinger equation with free Laplacian in an exterior domain. Our goal in Section~\ref{sec:equivalence} is to show that~\eqref{eq:classicalsmoothing} on the classical side and~\eqref{eq:Doi} on the quantum side can be deduced from each other through Egorov's theorem.

\begin{rmrk}[see {\cite[Remark 2.7]{P:23}}] \label{rem:classicalsmoothing}
Proposition~\ref{prop:classicalsmoothing} is a consequence of the more elementary estimate\cite[Corollary 2.6]{P:23}:
\begin{equation} \label{eq:classicalestimate}
\forall E \ge E_0, \forall r \ge 0, \forall (x, \xi) \in \{p = E\} , \qquad
	\textrm{Leb} \left( \left\{ t \in [0, T] \, : \, (\pi \circ \phi^t)(x, \xi) \in B_r(0) \right\} \right)
		\le C' \dfrac{r}{\sqrt{E}} .
\end{equation}
To understand the scaling of the variables $r$ and $E$ in~\eqref{eq:classicalestimate}, one can argue as follows: for fixed $r$, writing $\phi^t(x, \xi) = (x^t, \xi^t)$, we have in the high-energy limit $E \gg 1$:
\begin{equation*}
(\pi \circ \phi^t)(x, \xi) = x^t \in B_r(0)
	\qquad \Longrightarrow \qquad
E
	= p(x, \xi)
	= p\left(\phi^t(x, \xi)\right)
	= \tfrac{1}{2} |\xi^t|^2 + V(x^t)
	= \tfrac{1}{2} |\xi^t|^2 + O(1) ,
\end{equation*}
which means that the velocity of the trajectory is of order $|\xi^t| \approx \sqrt{2 E}$. Therefore, the curve $t \mapsto x^t$ can stay in $B_r(0)$ for a time of order at most $r/\sqrt{E}$. Going from~\eqref{eq:classicalestimate} to~\eqref{eq:classicalsmoothing} is a simple consequence of Fubini's theorem, seeing the radial map $x \mapsto \langle x \rangle^{- 2 \nu}$ as a weighted superposition of indicator functions of balls; see~\cite[Subsection 1.6]{P:23}.
\end{rmrk}

\subsection{Main result}
The proof of the equivalence of Theorem~\ref{thm:Doi} and Proposition~\ref{prop:classicalsmoothing} relies on Egorov's theorem. Introducing a semiclassical parameter $R \gg 1$, we obtain as a by-product estimates which relate the constants involved in $R$-dependent versions of the smoothing inequality~\eqref{eq:Doi} and of the escape rate estimate of the flow~\eqref{eq:classicalsmoothing}. Thus for any $R \ge 1$, we introduce two constants which quantify the classical escape rate estimate on the flow and the smoothing inequality respectively:
\begin{align}
\mathsf{C}_0(R)
	&:= \sup_{(x, \xi) \in \xR^{2d}} \sqrt{R^2 + p(x, \xi)} \int_0^T \dfrac{\xdif t}{\langle (\frac{1}{R} \pi \circ \phi^t) (x, \xi) \rangle^{2 \nu}} , \label{eq:classicalconstant}\\
\mathfrak{C}_0(R)
	&:= \sup_{u \in \xLtwo(\xR^d) \setminus \{0\}} \dfrac{1}{\|u\|_{\xLtwo(\xR^d)}^2} \int_0^T \left\| \langle x/R \rangle^{-\nu} \mathrm{Op}\left((R^2 + p)^{1/4}\right) e^{- i t P} u \right\|_{\xLtwo(\xR^d)}^2 \xdif t . \label{eq:quantumconstant}
\end{align}
Our main result reads as follows.
\begin{prpstn} \label{prop:constants}
Fix $T > 0$ and $\nu > 1/2$, and let $V$ satisfy Assumption~\ref{assum:V}. Then there exists a constant $c = c(T, V, \nu)$ such that the following inequalities hold:
\begin{equation} \label{eq:maininequalities}
\mathsf{C}_0(R) \left(1 + \dfrac{c}{R}\right)^{-1}
	\le \mathfrak{C}_0(R)
	\le \mathsf{C}_0(R) \left(1 + \dfrac{c}{R}\right) .
\end{equation}
\end{prpstn}
\begin{rmrk}
Notice that the quantity $\mathsf{C}_0(R)$ is non-decreasing in $R$. Estimating $\mathfrak{C}_0(R)$ for fixed $R$ (say $R = 1$), without considering an asymptotic regime such as $R \to + \infty$, seems more delicate and could possibly involve a combination of a quantitative analysis of high energies, using for instance the quantity
\begin{equation*}
\mathsf{C}_0^\infty
	:= \limsup_{(x, \xi) \to \infty} \sqrt{1 + p(x, \xi)} \int_0^T \dfrac{\xdif t}{\langle (\pi \circ \phi^t) (x, \xi) \rangle^{2 \nu}} ,
\end{equation*}
and low energies, using quantitative unique continuation. We refer to the paper of Laurent and Léautaud \cite{LL:16} for a study of related questions concerning observability cost and optimal observation time in the setting of the observability of the wave equation in compact Riemannian manifolds. See also the discussion in~\cite[Subsection 1.2]{P:23} in the context of the observability of the Schrödinger equation in the Euclidean space.
\end{rmrk}

\subsection{Smoothing effect and observability}

As already mentioned in~\cite[Subsection 1.6]{P:23}, smoothing inequalities associated with some operator $P$ are related to understanding which subsets of $\xR^d$ fail to observe the Schrödinger equation, in a quantitative way. We recall that a (measurable) subset $\omega \subset \xR^d$ observes the Schrödinger equation~\eqref{eq:Sch} in time $T > 0$ if
\begin{equation} \label{eq:obs}
\exists C > 0 : \forall u \in \xLtwo(\xR^d) , \qquad
	\left\| u \right\|_{\xLtwo(\xR^d)}^2
		\le C \int_0^T \left\| e^{- i t P} u \right\|_{\xLtwo(\omega)}^2 \xdif t .
\end{equation}
The meaning of this observability inequality is that the set $\omega$ captures a fraction of the energy of the solution $t \mapsto e^{- i t P} u$ over the time interval $[0, T]$, uniformly with respect to the initial state $u$. This inequality is somewhat similar to the smoothing inequality~\eqref{eq:Doi}: the indicator function $\mathbf{1}_\omega(x)$, seen as a function on phase space, plays the role of $p(x, \xi)^{1/4} \langle x \rangle^{- \nu}$. The major difference of course is that the smoothing inequality~\eqref{eq:Doi} goes the other way around compared to the observability inequality~\eqref{eq:obs}.

\section{Quantum-classical correspondence principle and proof of Proposition~\ref{prop:constants}} \label{sec:equivalence}

The main ingredient of the proof of Proposition~\ref{prop:constants} is a slightly non standard version of Egorov's theorem, a result at the core of the quantum-classical correspondence principle. See~\cite[Chapter 11]{Zworski} for a general introduction to Egorov's theorem in the semiclassical setting, or~\cite[Section 6]{Taylor1}, \cite[Chapter 7, Section 8]{Taylor2} for a microlocal version. See also~\cite{P:24Egorovinprep} and references therein.

\subsection{Quantization} \label{subsec:quantization}

Microlocal analysis relies on the concept of quantization, that is to say a correspondence between (smooth) functions on phase space (classical observables or symbols), and operators acting on $\xLtwo(\xR^d)$ (quantum observables).
Here we work with the Weyl quantization: for any function $a$ in the Schwartz space $\mathcal{S}(\xR^{2d})$, the Weyl quantization of $a$ is defined as the operator
\begin{equation*}
\left[\mathrm{Op}\left(a\right) u\right](x)
	= (2 \pi)^{-d} \int_{\xR^d} \int_{\xR^d} e^{i \xi \cdot (x - y)} a\left(\dfrac{x + y}{2}, \xi\right) u(y) \xdif y \xdif \xi ,
		\qquad x \in \xR^d , u \in \mathcal{S}(\xR^d) .
\end{equation*}
Such operators are called pseudo-differential operators. Classical references are~\cite{Hormander3,Martinez,Lerner} in the microlocal setting, or~\cite{Zworski,DimassiSjostrand} in the semiclassical setting. The quantization procedure $a \mapsto \mathrm{Op}(a)$ is linear and can be extended to $a \in \mathcal{S}'(\xR^{2d})$. Intuitively, it behaves as follows: classical observables $a$ depending only on the $x$ variable are mapped to the operator of multiplication by $a$, while symbols depending only on the $\xi$ variable are mapped to the corresponding Fourier multiplier. For general $a$ depending on both position and momentum variables, the operator $\mathrm{Op}(a)$ interpolates as much as possible between the two previous situations. Notice that we have in particular $P = \mathrm{Op}(p)$.

Symbols are generally considered as elements of a symbol class, a Fréchet space that collects symbols with specific growth or decay properties.

\begin{dfntn}[Order function] \label{def:orderfunction}
We say a function $f : \xR^{2d} \to \xR_+^\star$ is an order function, if it satisfies
\begin{equation} \label{eq:orderfunction}
\exists C > 0, N \ge 0 : \forall (x_0, \xi_0), (x, \xi) \in \xR^{2d}, \qquad
	f(x_0 + x, \xi_0 + \xi)
		\le C f(x_0, \xi_0) \langle (x, \xi) \rangle^N .
\end{equation}
The constants $C$ and $N$ are called the structure constants of $f$.
\end{dfntn}

The following symbol classes are not standard in semiclassical or microlocal analysis, but they fit in the Weyl--Hörmander calculus framework~\cite{Lerner}.

\begin{dfntn}[Symbol classes]
Let $f : \xR^{2d} \to \xR_+^\star$ be an order function.
We say a smooth function $a : \xR^{2d} \to \xC$ belongs to the class $S(f)$ if
\begin{equation*}
\forall \alpha \in \xN^{2d}, \exists C_\alpha > 0 : \forall (x, \xi) \in \xR^{2d}, \qquad
	\left| \partial^\alpha  a(x, \xi) \right|
		\le C_\alpha f(x, \xi) .
\end{equation*}
The space $S(f)$ is equipped with a Fréchet space structure with the seminorms
\begin{equation*}
|a|_{S(f)}^{(\ell)}
	= \max_{\substack{\alpha \in \xN^{2d} \\ 0 \le |\alpha| \le \ell}} \sup_{(x, \xi) \in \xR^{2d}} \left| \dfrac{\partial^\alpha a(x, \xi)}{f(x, \xi)} \right| .
\end{equation*}
\end{dfntn}

\begin{rmrk}
In the sequel, we will write $a \in S(f)$ for a vector-valued or matrix-valued function $a$ to mean that all its components belong to $S(f)$.
\end{rmrk}

Classical results of microlocal analysis are recalled in the Appendix~\ref{app}. They concern the composition of pseudo-differential operators (i.i.\ pseudo-differentiential calculus, see Proposition~\ref{prop:pseudocalc}), the $\xLtwo$ boundedness properties of these operators (Calder\'{o}n--Vaillancourt theorem, see Proposition~\ref{prop:CV}), and the sharp G{\aa}rding inequality (Proposition~\ref{prop:Gaarding}).

\subsection{Egorov's theorem}

The key argument of this note consists in applying a suitable version of Egorov's theorem, which roughly asserts that
\begin{equation} \label{eq:ideaEgorov}
e^{i t P} \mathrm{Op}\left(a\right) e^{- i t P}
	\approx \mathrm{Op}\left(a \circ \phi^t\right) .
\end{equation}
The version of Egorov's theorem that we use is deduced from~\cite{P:24Egorovinprep}.

\begin{thrm} \label{prop:Egorov}
Suppose $V$ is subject to Assumption~\ref{assum:V}. Fix $T > 0$ and let $f : \xR^{2d} \to \xR_+^\star$ be an order function. Assume that $f \circ \phi^t$ is an order function with uniform structure constants (see Definition~\ref{def:orderfunction}) with respect to $t \in [0, T]$. Then the following holds: for all $a \in S(f)$, we have
\begin{equation*}
e^{i t P} \mathrm{Op}\left(a\right) e^{- i t P}
	= \mathrm{Op}\left(a \circ \phi^t + \mathcal{R}_a(t)\right) ,
		\qquad t \in [0, T] ,
\end{equation*}
where $a \circ \phi^t \in S(f \circ \phi^t)$, $\mathcal{R}_a(t) \in S(f \circ \phi^t)$ for all $t \in [0, T]$ and
\begin{align}
\forall \ell \in \xN , \exists C_\ell > 0, \exists k \in \xN : \forall t \in [0, T] , &&
	\left| a \circ \phi^t \right|_{S(f \circ \phi^t)}^{(\ell)}
		&\le C_\ell \left| a \right|_{S(f)}^{(k)} , \label{eq:remainderboundedegorov}\\
\forall \ell \in \xN , \exists C_\ell > 0, \exists k \in \xN : \forall t \in [0, T] , &&
	\left| \mathcal{R}_a(t) \right|_{S(f \circ \phi^t)}^{(\ell)}
		&\le C_\ell \left| \nabla a \right|_{S(f)}^{(k)} .\label{eq:remainderboundedegorovR}
\end{align}
The constants $C_\ell$ depend on $f$ only through structure constants (Definition~\ref{def:orderfunction}).
\end{thrm}

This result relies on a general version of Egorov's theorem in the Weyl--Hörmander framework of metrics on the phase space, whose proof is quite involved. In fact, the key point of Theorem~\ref{prop:Egorov} is that the conjugated operator $e^{i t P} \mathrm{Op}(a) e^{- i t P}$ is a pseudo-differential operator, whose symbol lies in a symbol class associated with the initial order function composed by the flow. We now briefly give a proof of Theorem~\ref{prop:Egorov} from~\cite{P:24Egorovinprep}, using the notation in that reference.

\begin{proof}[Proof of Theorem~\ref{prop:Egorov} from~{\cite{P:24Egorovinprep}}]
We apply~\cite[Proposition 1.25]{P:24Egorovinprep} with $\hbar = 1$, $\beta = 0$ and the potential $V$ which is sub-quadratic in the sense that $\nabla^k V$ is bounded for all $k \ge 2$. This shows that the assumptions of~\cite[Theorem 1.15]{P:24Egorovinprep} are verified, so that we have a version of Egorov's theorem in the symbol class $S(f) = S(f, g_{\mathrm{eucl}})$ (with $g_{\mathrm{eucl}} = \xdif x^2 + \xdif \xi^2$ the standard Euclidean metric), in a time window $[0, T]$ with fixed $T > 0$.
Then~\cite[Theorem 1.15]{P:24Egorovinprep} states that
\begin{equation*}
e^{i t P} \mathrm{Op}(a) e^{- i t P}
	= \mathrm{Op}(a \circ \phi^t) + \mathrm{Op}\left(\mathcal{R}_a(t)\right) ,
\end{equation*}
with the estimates~\eqref{eq:remainderboundedegorov} and~\eqref{eq:remainderboundedegorovR} provided $f \circ \phi^t$ is an admissible weight with respect to the metric $g_{\mathrm{eucl}}$. This is equivalent to saying that $f \circ \phi^t$ is an order function with uniform structure constants in the time interval $t \in [0, T]$, which we assume in the statement of Theorem~\ref{prop:Egorov}. This finishes the proof.
\end{proof}

\subsection{Technical results} \label{subsec:proofEgorov}

%We first prove a Lipschitz property of the flow.
%
%\begin{lmm} \label{lem:GronwallLip}
%Let $V$ satisfy~\eqref{eq:assumV}. Then there exists a constant $C > 0$ such that
%\begin{equation*}
%\forall t \in \xR, \forall \rho_1, \rho_2 \in \xR^{2d} , \qquad
%	\left|\phi^t(\rho_2) - \phi^t(\rho_1)\right|
%		\le e^{C |t|} \left| \rho_2 - \rho_1 \right| .
%\end{equation*}
%\end{lmm}
%
%\begin{proof}
%Set $H_p(x, \xi) = (\xi, - \nabla V(x))^\top$. We use Hamilton's equations~\eqref{eq:Hamilton} to write
%\begin{equation*}
%\left|\phi^t(\rho_2) - \phi^t(\rho_1)\right|
%	\le 	\left|\rho_2 - \rho_1\right| + \int_0^t \left|H_p\left(\phi^s(\rho_2)\right) - H_p\left(\phi^s(\rho_1)\right)\right| \xdif s .
%\end{equation*}
%Next we remark that $\nabla H_p$ is bounded uniformly on the whole phase space $\xR^{2d}$ by assumption on $V$~\eqref{eq:assumV}. Therefore the mean value inequality yields
%\begin{equation*}
%\left|\phi^t(\rho_2) - \phi^t(\rho_1)\right|
%	\le 	\left|\rho_2 - \rho_1\right| + \int_0^t \sup_{\xR^{2d}} |\nabla H_p|  \left|\phi^s(\rho_2) - \phi^s(\rho_1)\right| \xdif s .
%\end{equation*}
%The result follows by Grönwall's inequality.
%\end{proof}

We collect several technical results about the Hamiltonian flow and the symbolic properties of powers of $R^2 + p$ with $R \ge 1$.
We first recall known estimates on the derivatives of the Hamiltonian flow.

\begin{lmm}[Lemma 2.10 in {\cite{P:23}}] \label{lem:behaviorflow}
Suppose $V$ is subject to Assumption~\ref{assum:V} and fix $T > 0$. Then for all $k \in \xN \setminus \{0\}$, the {$k$-th} differential $\xdif^k \phi^t$ satisfies:
\begin{equation*}
\sup_{t \in [0, T]} \; \sup_{(x, \xi) \in \xR^{2d}} \left\| \xdif^k \phi^t (x, \xi) \right\| < \infty .
\end{equation*}
In particular, the map $\phi^t : \xR^{2d} \to \xR^{2d}$ is Lipschitz with a constant uniform in $t \in [0, T]$. In addition, the following holds:
\begin{equation*}
\forall T > 0, \forall \alpha \in \xN^{2d} \setminus \{0\} , \exists C_{\alpha, T} > 0 : \forall a \in S(1), \forall t \in [-T, T] , \qquad
	\| \partial^\alpha (a \circ \phi^t) \|_{\infty}
		\le C_{\alpha, T} \max_{\substack{\beta \in \xN^{2d} \\ 1 \le |\beta| \le |\alpha|}} \| \partial^\beta a \|_{\infty} .
\end{equation*}
\end{lmm}

In the sequel, we will consider a specific order function, namely
\begin{equation} \label{eq:f0}
f_R
	:= \langle x/R \rangle^{-2 \nu} \sqrt{R^2 + p} ,
		\qquad \nu > 1/2 , \, R \ge 1 .
\end{equation}
We observe that it blows up at infinity in $\xi$, which is why we need the refined version of Egorov's theorem given in~\cite{P:24Egorovinprep}. In addition, the extra parameter $R \ge 1$ plays the role of a semiclassical parameter and allows us to gain smallness in the remainder of Egorov's theorem. Let us prove that $f_R$ above satisfies the assumptions of Theorem~\ref{prop:Egorov}, uniformly with respect to $R \ge 1$.

\begin{lmm} \label{lem:prooff0}
Suppose $V$ is subject to Assumption~\ref{assum:V}. Then for any $\gamma \in \xR$, the map $(R^2 + p)^\gamma$ is an order function in the sense of Definition~\ref{def:orderfunction}. The map $f_R \circ \phi^t$ with $f_R$ defined in~\eqref{eq:f0} is an order function. Structure constants are uniform with respect to $t$ in any compact time interval and $R \ge 1$.
\end{lmm}

\begin{proof}
We denote phase space points by $\rho = (x, \xi) \in \xR^{2d}$, and we set $\tilde f(\rho) = \langle \pi(\rho)\rangle^{-2 \nu}$. 
Recall that $\langle x \rangle^{2m} + \langle \xi \rangle^2$ is an order function~\cite[4.4.1. Examples]{Zworski}. Since one can check from~\eqref{eq:assumV} that
\begin{equation*}
\dfrac{1}{C} \left(\langle x \rangle^{2m} + \langle \xi \rangle^2\right)
	\le 1 + p
	\le C \left(\langle x \rangle^{2m} + \langle \xi \rangle^2\right) ,
\end{equation*}
for some large enough constant $C > 0$, we deduce that $R^2 + p$ is an order function with uniform structure constants with respect to $R \ge 1$. The result for $(R^2 + p)^\gamma$ follows.

Now we check that $f_R \circ \phi^t$ is an order function. We use the fact that $\langle a + b \rangle \le 2 \langle a \rangle \langle b \rangle$ together with the Lipschitz property of the flow (Lemma~\ref{lem:behaviorflow}) to have
\begin{align*}
\left\langle \frac{1}{R} (\pi \circ \phi^t)(\rho_2) \right\rangle
	&= \left\langle \frac{1}{R} (\pi \circ \phi^t)(\rho_1) + \frac{1}{R} (\pi \circ \phi^t)(\rho_2) - \frac{1}{R} (\pi \circ \phi^t)(\rho_1) \right\rangle \\
	&\le 2 \left\langle \frac{1}{R} (\pi \circ \phi^t)(\rho_1) \right\rangle \left\langle \phi^t(\rho_2) - \phi^t(\rho_1) \right\rangle
	\le 2 C_T \left\langle \frac{1}{R} (\pi \circ \phi^t)(\rho_1) \right\rangle \left\langle \rho_2 - \rho_1 \right\rangle ,
\end{align*}
with a constant $C_T > 0$ uniform in $t \in [-T, T]$.
Taking the power $-2 \nu$, we obtain for any $\rho_1, \rho_2 \in \xR^{2d}$:
\begin{equation} \label{eq:tildef}
\left\langle \frac{1}{R} (\pi \circ \phi^t)(\rho_1) \right\rangle^{-2 \nu}
	\le (2 C_T)^{2 \nu} \left\langle \frac{1}{R} (\pi \circ \phi^t)(\rho_2) \right\rangle^{-2\nu} \left\langle \rho_2 - \rho_1 \right\rangle^{2 \nu} .
\end{equation}
Multiplying~\eqref{eq:tildef} by $(R^2 + p)^{1/2}$, which is an order function invariant by the flow~\eqref{eq:conservationofenergy}, it follows that $f_R \circ \phi^t$ is an order function with uniform structure constants with respect to $t \in [-T, T]$ and $R \ge 1$.
\end{proof}

Next we give some properties of the symbols relevant to this problem.

\begin{lmm} \label{lem:symbolclassesp}
Suppose $V$ is subject to Assumption~\ref{assum:V}. For all $\gamma \in (- \infty, 1/2]$, we have
\begin{equation} \label{eq:nablaqsymbolclass}
(R^2 + p)^\gamma \in S\left((R^2 + p)^\gamma\right)
	\qquad \textrm{and} \qquad
\nabla (R^2 + p)^\gamma \in S\left((R^2 + p)^{\gamma - 1/2}\right) \subset S\left(\frac{1}{R} (R^2 + p)^\gamma\right) ,
\end{equation}
uniformly with respect to $R \ge 1$.
\end{lmm}

\begin{proof}
It is sufficient to check the desired estimate~\eqref{eq:nablaqsymbolclass} on $\nabla (R^2 + p)^\gamma$. We see $(R^2 + p)^\gamma$ as the composition of $s \mapsto (R^2 + s)^\gamma$ with $p$. From the Faà di Bruno formula, for any non-zero $\alpha \in \xN^{2d}$, the derivative $\partial^\alpha (R^2 + p)^\gamma$ is a sum of terms of the form
\begin{equation} \label{eq:termsinfaadibruno}
(R^2 + p)^{\gamma - j} \prod_{k = 1}^j \partial^{\beta_k} p ,
	\qquad \textrm{with} \;\, j \in \{1, 2, \ldots, |\alpha|\} \;\, \textrm{and} \;\, \sum_{k=1}^j \beta_k = \alpha , \;\, \beta_k \in \xN^{2d} \setminus \{0\}, \, \forall k ,
\end{equation}
up to combinatorial constants. In view of Assumption~\ref{assum:V} on the potential $V$, we have $|\partial^\beta p| \le C_\beta (1 + p)^{1/2}$ for any $\beta \in \xN^{2d} \setminus \{0\}$ (notice that $m \le 1$ in~\eqref{eq:assumV} is important here), so that terms of the form~\eqref{eq:termsinfaadibruno} are bounded by
\begin{equation} \label{eq:estimatederivatives}
(R^2 + p)^{\gamma - j} (1 + p)^{j/2}
	= (R^2 + p)^{\gamma} \left( \dfrac{\sqrt{1 + p}}{R^2 + p} \right)^j
	\le (R^2 + p)^{-(1/2 - \gamma)}
	\le \dfrac{1}{R} (R^2 + p)^\gamma
\end{equation}
(we used that $j \ge 1$ from~\eqref{eq:termsinfaadibruno} in the first inequality). Therefore $\nabla (R^2 + p)^\gamma$ belongs to the desired symbol class~\eqref{eq:nablaqsymbolclass}.
\end{proof}

\subsection{Reduction step}

%\begin{rmrk} \label{rem:seminormsvanish}
%Examining the details of the proof, and more specifically the last inequality of~\eqref{eq:estimatederivatives}, one observes that the seminorms of $\nabla (c + p)^\gamma$ in the symbol class $S(1)$ are in fact bounded from above by $c^{\gamma - 1/2}$, which tends to zero for $\gamma < 1/2$. This will be important in Lemma~\ref{lem:positivity} below.
%\end{rmrk}

The norm squared in the left-hand side of~\eqref{eq:Doi} can be rewritten as
\begin{equation*}
\int_0^T \left(u, e^{i t P} \mathrm{Op}\left((1 + p)^{1/4}\right) \langle x \rangle^{-2 \nu} \mathrm{Op}\left((1 + p)^{1/4}\right) e^{- i t P} u \right)_{\xLtwo} \xdif t ,
\end{equation*}
where $(\bullet, \bullet)_{\xLtwo}$ is the $\xLtwo$ inner product.
Our first step is to rewrite the operator in the middle as a genuine pseudo-differential operator, up to a bounded operator. Recall that we are concerned with an $R$-dependent version of~\eqref{eq:Doi} (see~\eqref{eq:quantumconstant}). So from now on, we write for convenience
\begin{equation} \label{eq:defQ}
Q = Q_R
	:= \mathrm{Op}\left((R^2 + p)^{1/4}\right) .
\end{equation}

\begin{lmm} \label{lem:Q}
Suppose $V$ is subject to Assumption~\ref{assum:V}. For any $R \ge 1$, there exists a bounded self-adjoint operator $A_R$ such that
\begin{equation} \label{eq:pseudocalcandAR}
Q \langle x/R \rangle^{-2 \nu} Q
	= \mathrm{Op}\left(\dfrac{\sqrt{R^2 + p}}{\langle x/R \rangle^{2 \nu}}\right) + \dfrac{1}{R} A_R ,
\end{equation}
as operators acting on $\mathcal{S}(\xR^d)$, with $A_R$ bounded uniformly with respect to $R \ge 1$.
\end{lmm}

\begin{proof}
%In this proof, we do not care about the dependence on the parameter $R$ (in particular the magnitude of $\|A_R\|_{\xLtwo \to \xLtwo}$ is irrelevant).
We first apply the pseudo-differential calculus (Proposition~\ref{prop:pseudocalc}~\eqref{eq:pseudocalc1}) to have
\begin{equation} \label{eq:1ststeppseudocalc}
Q \langle x/R \rangle^{-2 \nu}
	= \mathrm{Op}\left(\dfrac{(R^2 + p)^{1/4}}{\langle x/R \rangle^{2 \nu}}\right) + \mathrm{Op}\left(r_1\right) ,
\end{equation}
where $r_1 \in S(\frac{1}{R} (R^2 + p)^{-1/4})$, since Lemma~\ref{lem:symbolclassesp} implies
\begin{equation*}
\nabla (R^2 + p)^{1/4} \in S((R^2 + p)^{-1/4})
	\qquad \textrm{and} \qquad
\nabla \langle x/R \rangle^{-2 \nu} \in S(\frac{1}{R}) \,.
\end{equation*}
Now we wish to compose~\eqref{eq:1ststeppseudocalc} on the right with $Q$. We have
\begin{equation*}
 \mathrm{Op}\left(\dfrac{(R^2 + p)^{1/4}}{\langle x/R \rangle^{2 \nu}}\right) Q
 	= \mathrm{Op}\left(\dfrac{(R^2 + p)^{1/2}}{\langle x/R \rangle^{2 \nu}}\right) + \mathrm{Op}(r_2) ,
\end{equation*}
where $r_2 \in S(1/R)$, uniformly with respect to $R \ge 1$, since on the one hand
\begin{equation*}
\nabla (R^2 + p)^{1/4} \in S((R^2 + p)^{-1/4}) ,
\end{equation*}
in virtue of Lemma~\ref{lem:symbolclassesp}, and on the other hand
\begin{align} \label{eq:symbolclassnablaf_R}
\nabla \dfrac{(R^2 + p)^{1/4}}{\langle x/R \rangle^{2 \nu}}
	&= \langle x/R \rangle^{-2 \nu} \nabla (R^2 + p)^{1/4} + (R^2 + p)^{1/4} \times \dfrac{1}{R} \left(\nabla \langle \bullet \rangle^{-2\nu}\right)(x/R) \nonumber\\
	&\in S\left(\dfrac{1}{R} (R^2 + p)^{1/4}\right) .
\end{align}
Another application of Proposition~\ref{prop:pseudocalc}~\eqref{eq:pseudocalc1} yields
\begin{equation*}
\mathrm{Op}(r_1) Q
	= \mathrm{Op}(b)
\end{equation*}
with $b \in S(1/R)$ since
\begin{equation*}
r_1 \in S(\frac{1}{R} (R^2 + p)^{-1/4})
	\qquad \textrm{and} \qquad
(R^2 + p)^{1/4} \in S((R^2 + p)^{1/4}) .
\end{equation*}
In conclusion, we have
\begin{equation*}
Q \langle x/R \rangle^{-2 \nu} Q
	= \mathrm{Op}\left(\dfrac{(R^2 + p)^{1/2}}{\langle x/R \rangle^{2 \nu}} + r_2 + b\right) ,
\end{equation*}
where $r_2 + b \in S(1/R)$. The Calder{\'o}n--Vaillancourt theorem (Proposition~\ref{prop:CV}) gives the result.
\end{proof}

We now apply Egorov's theorem (Theorem~\ref{prop:Egorov}). The (large) parameter $R \ge 1$ will be useful to justify that the remainder in Egorov's theorem (which is a priori an unbounded operator) is in some sense smaller than the leading term.

\begin{lmm} \label{lem:key}
Fix $T > 0$ and $\nu > 1/2$. Then for any $R \ge 1$ and for all $t \in [0, T]$, there exists a bounded self-adjoint operator $A_R(t)$ (bounded uniformly in $t \in [0, T]$ and $R \ge 1$) and $c_R(t) \in S(f_R \circ \phi^t)$ such that
\begin{equation} \label{eq:inequalities}
e^{i t P} Q \langle x/R \rangle^{-2 \nu} Q e^{- i t P}
	= \mathrm{Op}\left(f_R \circ \phi^t\right) + \dfrac{1}{R} \mathrm{Op}\left(c_R(t)\right) + \dfrac{1}{R} A_R(t)  ,
\end{equation}
where
\begin{equation*}
\forall \ell \in \xN, \qquad
	\sup_{\substack{R \ge 1 \\ t \in [0, T]}} \left|c_R(t)\right|_{S(f_R \circ \phi^t)}^{(\ell)} < \infty .
\end{equation*}
\end{lmm}

\begin{proof}
We conjugate~\eqref{eq:pseudocalcandAR} by the Schrödinger propagator to obtain:
\begin{equation} \label{eq:intermediatestep}
e^{i t P} Q \langle x/R \rangle^{-2 \nu} Q e^{- i t P}
	= e^{i t P} \mathrm{Op}\left(f_R\right) e^{- i t P} + \dfrac{1}{R} e^{i t P} A_R e^{-i t P} .
\end{equation}

To prove~\eqref{eq:inequalities}, it remains to apply Egorov's theorem (Theorem~\ref{prop:Egorov}) to the first term in the right-hand side of~\eqref{eq:intermediatestep}.
Theorem~\ref{prop:Egorov} applies here with $a = f_R \in S(f_R)$. Notice that $f_R \circ \phi^t$ is an order function with structure constants uniform in $t \in [0, T]$ and $R \ge 1$ in virtue of Lemma~\ref{lem:prooff0}. We have
\begin{equation} \label{eq:egorovonfR}
e^{i t P} \mathrm{Op}\left(f_R\right) e^{- i t P}
	= \mathrm{Op}\left(f_R \circ \phi^t + \mathcal{R}(t)\right)
\end{equation}
where $f_R \circ \phi^t$ lies in a bounded subset of $S(f_R \circ \phi^t)$ by~\eqref{eq:remainderboundedegorov} and
\begin{equation} \label{eq:estremainderR(t)}
\sup_{t \in [0, T]} \left|\mathcal{R}(t)\right|_{S(f_R \circ \phi^t)}^{(\ell)}
	\le C_\ell \left|\nabla f_R\right|_{S(f_R)}^{(\ell)} .
\end{equation}
The computation~\eqref{eq:symbolclassnablaf_R} shows that $\nabla f_R \in \frac{1}{R} S(f_R)$, uniformly with respect to $R \ge 1$ (a consequence of Lemma~\ref{lem:symbolclassesp} applied with $\gamma = 1/2$). This proves ultimately that $\mathcal{R}(t) \in \frac{1}{R} S(f_R \circ \phi^t)$, uniformly with respect to $t \in [0, T]$ and to the parameter $R \ge 1$.

Combining~\eqref{eq:intermediatestep} and \eqref{eq:egorovonfR}, we obtain
\begin{equation*}
e^{i t P} Q \langle x/R \rangle^{-2 \nu} Q e^{- i t P}
	= \mathrm{Op}\left(f_R \circ \phi^t\right) + \mathrm{Op}\left(\mathcal{R}(t)\right) + \dfrac{1}{R} e^{i t P} A_R e^{-i t P} .
\end{equation*}
We obtain the thought equality~\eqref{eq:inequalities} by setting $A_R(t) = e^{i t P} A_R e^{-i t P}$, which is uniformly bounded in $t \in \xR$ and $R \ge 1$, and $c_R(t) = R \mathcal{R}(t)$, which belongs to $S(f_R \circ \phi^t)$ uniformly in $R \ge 1$ and $t \in [0, T]$.
\end{proof}

\subsection{Proof of Proposition~\ref{prop:constants}}

With Lemma~\ref{lem:key} at hand, we can proceed with the proof Proposition~\ref{prop:constants}, namely the quantitative version of the equivalence of Theorem~\ref{thm:Doi} and Proposition~\ref{prop:classicalsmoothing}. We set
\begin{equation} \label{eq:thesymbol}
a_R : (x, \xi) \longmapsto \int_0^T \dfrac{(R^2 + p)^{1/2}}{\langle \frac{1}{R} \pi \circ \phi^t \rangle^{2 \nu}} \xdif t
	= \int_0^T f_R \circ \phi^t \xdif t ,
\end{equation}
where $f_R$ is defined in~\eqref{eq:f0}.

\begin{proof}[Proof the right-hand side of~\eqref{eq:maininequalities}]
By definition, the symbol $a_R$ defined in~\eqref{eq:thesymbol} is bounded by the constant $\mathsf{C}_0(R)$, defined in~\eqref{eq:classicalconstant}. Our goal is to deduce the smoothing inequality with parameter $R$ (see~\eqref{eq:quantumconstant}) through Egorov's theorem.
In view of Lemma~\ref{lem:key}, we have
\begin{multline*}
\int_0^T e^{i t P} \mathrm{Op}\left((R^2 + p)^{1/4} \langle x/R \rangle^{-2 \nu} (R^2 + p)^{1/4}\right) e^{- i t P} \xdif t \\
	= \mathrm{Op}\left(\int_0^T \left(f_R \circ \phi^t + \frac{1}{R} c_R(t)\right) \xdif t\right) + \int_0^T \frac{1}{R} A_R(t) \xdif t ,
\end{multline*}
with $A_R(t)$ bounded uniformly with respect to $t \in [0, T]$ and $R \ge 1$ and $\int_0^T c_R(t) \xdif t \in S(a_R)$ (uniformly in $R \ge 1$).
Thus, all we need is to check that the symbol $a_R = \int_0^T f_R \circ \phi^t \xdif t$ belongs to the symbol class $S(\mathsf{C}_0(R))$ as well (uniformly in $R \ge 1$), and the smoothing inequality will follow from the sharp G{\aa}rding inequality (Proposition~\ref{prop:Gaarding}).
We know from Egorov's theorem (Theorem~\ref{prop:Egorov}) that the seminorms of $f_R \circ \phi^t$ in $S(f_R \circ \phi^t)$ are bounded uniformly with respect to $t \in [0, T]$ and $R \ge 1$.
Integrating in time, we deduce that for all $(x, \xi) \in \xR^{2d}$ and all $R \ge 1$,
\begin{equation*}
\left| \partial^\alpha \int_0^T (f_R \circ \phi^t) (x, \xi) \xdif t \right|
	\le C_{\alpha, T} \int_0^T (f_R \circ \phi^t) (x, \xi) \xdif t
	\le C_{\alpha, T} \mathsf{C}_0(R) ,
\end{equation*}
by differentiating under the inegral sign. This means that the symbol $a_R$ defined in~\eqref{eq:thesymbol} is indeed in $S(\mathsf{C}_0(R))$ uniformly with respect to $R \ge 1$. From~\eqref{eq:symbolclassnablaf_R}, we recall that $\nabla f_R \in S(\frac{1}{R} f_R)$ uniformly with respect to $R \ge 1$, so we deduce from Lemma~\ref{lem:behaviorflow} that
\begin{equation} \label{eq:gainfornabla}
\nabla (f_R \circ \phi^t) \in S\left(\tfrac{1}{R} f_R \circ \phi^t\right) ,
\end{equation}
uniformly with respect to $R \ge 1$ and $t \in [-T, T]$. Thus the sharp G{\aa}rding inequality (Proposition~\ref{prop:Gaarding}) applied to $\mathsf{C}_0(R) - a_R$ yields
\begin{equation*}
\mathrm{Op}(a_R)
	= \mathrm{Op}\left(\int_0^T f_R \circ \phi^t \xdif t\right)
 	\le \mathsf{C}_0(R) + \dfrac{c}{R} \mathsf{C}_0(R) ,
\end{equation*}
for some constant $c$ independent of $R$. Furthermore, the Calder\'{o}n--Vaillancourt theorem (Proposition~\ref{prop:CV}) applied to $\int_0^T c_R(t) \xdif t \in S(\mathsf{C}_0(R))$, allows us to conclude that
\begin{equation*}
\int_0^T e^{i t P} \mathrm{Op}\left((R^2 + p)^{1/4} \langle x/R \rangle^{-2 \nu} (R^2 + p)^{1/4}\right) e^{- i t P} \xdif t
	\le \mathsf{C}_0(R) + \dfrac{c'}{R} \left(1 + \mathsf{C}_0(R)\right) ,
\end{equation*}
hence $\mathfrak{C}_0(R) \le \mathsf{C}_0(R) (1 + O(\frac{1}{R}))$.
\end{proof}

\begin{proof}[Proof the left-hand side of~\eqref{eq:maininequalities}]
Our goal is to establish an estimate on the $R$-dependent classical averages of the flow in~\eqref{eq:classicalconstant} in terms of $\mathfrak{C}_0(R)$ defined in~\eqref{eq:quantumconstant}. We apply Lemma~\ref{lem:key} to rewrite the left-hand side of the smoothing inequality of Theorem~\ref{thm:Doi} as:
\begin{equation*}
\int_0^T \left( u, e^{i t P} Q \left\langle x/R \right\rangle^{-2 \nu} Q e^{-i t P} u \right)_{\xLtwo} \xdif t
	= \int_0^T \left( u, \mathrm{Op}\left(f_R \circ \phi^t + \dfrac{1}{R} c_R(t) + \dfrac{1}{R} A_R(t)\right) u \right)_{\xLtwo} \xdif t ,
\end{equation*}
for all $u \in \xLtwo(\xR^d)$, and we obtain by definition of~\eqref{eq:quantumconstant}:
\begin{equation} \label{eq:fnuRplusremainder}
\int_0^T \left( u, \mathrm{Op}\left(f_R \circ \phi^t + \dfrac{1}{R} c_R(t)\right) u \right)_{\xLtwo} \xdif t
	\le \left(\mathfrak{C}_0(R) + \dfrac{1}{R} \sup_{\substack{R \ge 1 \\ t \in [0, T]}} \left\| A_R(t) \right\|_{\xLtwo \to \xLtwo}\right) \|u\|_{\xLtwo(\xR^d)}^2 ,
\end{equation}
with $c_R(t) \in S(f_R \circ \phi^t)$ uniformly with respect to $R \ge 1$ and $t \in [0, T]$.

Now we apply the inequality~\eqref{eq:fnuRplusremainder} to a particular wave function, namely the Gaussian wave packet microlocalized near a phase space point $(x_0, \xi_0) \in \xR^{2d}$:
\begin{equation*}
u(x)
	= \pi^{-d/4} \exp\left(- \dfrac{|x - x_0|^2}{2}\right) e^{i \xi_0 \cdot x} ,
		\qquad x \in \xR^d .
\end{equation*}
It is normalized $\| u \|_{\xLtwo} = 1$ and a standard computation~\cite[Proposition (1.48)]{folland} shows that for any symbol $a \in S(f)$ (where $f$ is any order function), it holds
\begin{equation*}
\left( u, \mathrm{Op}\left(a\right) u \right)_{\xLtwo}
	= \pi^{-d} \int_{\xR^{2d}} a(x + x_0, \xi + \xi_0) \exp\left(- |(x, \xi)|^2\right) \xdif x \xdif \xi .
\end{equation*}
Taking $a := \int_0^T (f_R \circ \phi^t + \frac{1}{R} c_R(t)) \xdif t$, one deduces from~\eqref{eq:fnuRplusremainder} that for all $(x_0, \xi_0) \in \xR^{2d}$,
\begin{equation} \label{eq:C_0}
\int_{\xR^{2d}} a(x + x_0, \xi + \xi_0) \pi^{-d} \exp\left(- |(x, \xi)|^2\right) \xdif x \xdif \xi
	\le \mathfrak{C}_0(R) + \dfrac{1}{R} \sup_{\substack{R \ge 1 \\ t \in [0, T]}} \left\| A_R(t) \right\|_{\xLtwo \to \xLtwo} .
\end{equation}
Yet the left-hand side can be written as:
\begin{align}
\int_{\xR^{2d}} a(x + x_0, \xi + \xi_0) \pi^{-d} \exp\left(- |(x, \xi)|^2\right) \xdif x \xdif \xi
	&= \int_0^T (f_R \circ \phi^t)(x_0, \xi_0) \xdif t \label{eq:firstterm} \\
		&\hspace*{-4cm}+ \int_0^T \int_{\xR^{2d}} \int_0^1 \nabla (f_R \circ \phi^t)(x_0 + s x, \xi_0 + s \xi) \cdot (x, \xi) \pi^{-d} \exp\left(- |(x, \xi)|^2\right) \xdif s \xdif x \xdif \xi \xdif t \label{eq:secondterm} \\
		&\hspace*{-4cm}+ \dfrac{1}{R} \int_0^T \int_{\xR^{2d}} c_R(t)(x_0 + x, \xi_0 + \xi) \pi^{-d} \exp\left(- |(x, \xi)|^2\right) \xdif x \xdif \xi \xdif t . \label{eq:thirdterm}
\end{align}
The last term~\eqref{eq:thirdterm} in the right-hand side is controlled as follows: since $c_R(t) \in S(f_R \circ \phi^t)$ and $f_R \circ \phi^t$ is an order function with structure constants uniform in $t \in [0, T]$ and $R \ge 1$ (Lemma~\ref{lem:prooff0}), we deduce that
\begin{align} \label{eq:estthirdterm}
\left|\int_0^T \int_{\xR^{2d}} c_R(t)(x_0 + x, \xi_0 + \xi) \pi^{-d} \exp\left(- |(x, \xi)|^2\right) \xdif x \xdif \xi \xdif t\right| & \nonumber\\
	&\hspace*{-5cm}\le C \int_0^T \int_{\xR^{2d}} (f_R \circ \phi^t)(x_0 + x, \xi_0 + \xi) \pi^{-d} \exp\left(- |(x, \xi)|^2\right) \xdif x \xdif \xi \xdif t \nonumber\\
	&\hspace*{-5cm}\le C' \int_0^T \int_{\xR^{2d}} (f_R \circ \phi^t)(x_0, \xi_0) \pi^{-d} \langle (x, \xi)\rangle^N \exp\left(- |(x, \xi)|^2\right) \xdif x \xdif \xi \xdif t \nonumber\\
	&\hspace*{-5cm}\le C'' \int_0^T (f_R \circ \phi^t)(x_0, \xi_0) \xdif t ,
\end{align}
with a constant $C''$ independent of $R \ge 1$. Now we handle the second term~\eqref{eq:secondterm} as follows: we recall that $\nabla (f_R \circ \phi^t) \in S\left(\tfrac{1}{R} f_R \circ \phi^t\right)$ uniformly with respect to $R \ge 1$ and $t \in [-T, T]$ (see~\eqref{eq:gainfornabla}). Once again we use the fact that $f_R \circ \phi^t$ is an order function with structure constants uniform in $R \ge 1$ and $t \in [-T, T]$ to obtain
\begin{align} \label{eq:estsecondterm}
\left|\int_0^T \int_{\xR^{2d}} \int_0^1 \nabla (f_R \circ \phi^t)(x_0 + s x, \xi_0 + s \xi) \cdot (x, \xi) \pi^{-d} \exp\left(- |(x, \xi)|^2\right) \xdif s \xdif x \xdif \xi \xdif t\right| & \nonumber\\
	&\hspace*{-9cm}\le \dfrac{C}{R} \int_0^T \int_{\xR^{2d}} \int_0^1 (f_R \circ \phi^t)(x_0, \xi_0) |s|^N  \langle (x, \xi) \rangle^{1 + N} \pi^{-d} \exp\left(- |(x, \xi)|^2\right) \xdif s \xdif x \xdif \xi \xdif t \nonumber\\
	&\hspace*{-9cm}\le \dfrac{C'}{R} \int_0^T (f_R \circ \phi^t)(x_0, \xi_0) \xdif t .
\end{align}
Going back to~\eqref{eq:firstterm}, \eqref{eq:secondterm}, \eqref{eq:thirdterm}, we conclude that
\begin{equation*}
\int_{\xR^{2d}} a(x + x_0, \xi + \xi_0) \pi^{-d} \exp\left(- |(x, \xi)|^2\right) \xdif x \xdif \xi
	= \left( 1 + O\left(\dfrac{1}{R}\right) \right) \int_0^T (f_R \circ \phi^t)(x_0, \xi_0) \xdif t ,
\end{equation*}
and the constant involved in the $O$ notation is independent of $(x_0, \xi_0) \in \xR^{2d}$. Plugging this into~\eqref{eq:C_0}, we finally obtain
\begin{equation*}
\int_0^T (f_R \circ \phi^t)(x_0, \xi_0) \xdif t
	\le \dfrac{\mathfrak{C}_0(R) + O(\frac{1}{R})}{1 + O(\frac{1}{R})}
	= \mathfrak{C}_0(R) \left( 1 + O\left(\dfrac{1}{R}\right) \right) .
\end{equation*}
Therefore $\mathsf{C}_0(R) \le \mathfrak{C}_0(R) (1 + O(1/R))$, which concludes the proof of Proposition~\ref{prop:constants}.
\end{proof}

\subsection*{Acknowledgements} The author thanks Matthieu Léautaud for many helpful discussions, pointing references and reading a preliminary version of this article. The author is also grateful to the anonymous referee for his/her valuable comments on the first version, which led to a significant improvement of the article.

\appendix

\section{Toolbox of microlocal analysis} \label{app}

\begin{prpstn}[Pseudo-differential calculus, see {\cite[Proposition B.5]{P:23}}] \label{prop:pseudocalc}
Let $f_1, f_2$ be two order functions. Then for all $a_1, a_2 \in \xCinfty(\xR^{2d})$ such that $\nabla a_1 \in S(f_1)$ and $\nabla a_2 \in S(f_2)$, we have
\begin{equation*}
\mathrm{Op}\left(a_1\right) \mathrm{Op}\left(a_2\right)
	= \mathrm{Op}\left(a_1 a_2\right) + \mathrm{Op}\left(\mathcal{R}_1(a_1, a_2)\right) ,
\end{equation*}
where
\begin{equation} \label{eq:pseudocalc1} 
\forall \ell \in \xN, \exists k \in \xN, \exists C_\ell > 0 : \qquad
	\left| \mathcal{R}_1(a_1, a_2) \right|_{S(f_1 f_2)}^{(\ell)}
		\le C_\ell \left| \nabla a_1 \right|_{S(f_1)}^{(k)} \left| \nabla a_2 \right|_{S(f_2)}^{(k)} .
\end{equation}
\end{prpstn}

We refer to~\cite[Appendix B]{P:23} for a proof. Notice that the statement of~\cite[Proposition B.5]{P:23} is slightly different from the above since it works for symbols $(a_1, a_2) \in S(f_1) \times S(f_2)$. However, the proof remains valid for symbols $a_1$ and $a_2$ such that $\nabla a_j \in S(f_j)$, $j = 1, 2$. The operators $\mathrm{Op}(a_j)$ are properly defined since $a_j \in \mathcal{S}'(\xR^{2d})$; see~\cite[Definition 1.1.9]{Lerner}.

\begin{prpstn}[Calder\'{o}n--Vaillancourt, see {\cite{CV}}] \label{prop:CV}
There exist an integer $k$ and a constant $C > 0$ depending only on the dimension $d$ such that the following holds: for any $a \in S(1)$, the operator $\mathrm{Op}(a)$ extends to a bounded operator on $\xLtwo(\xR^d)$ with
\begin{equation} \label{eq:CV}
\left\| \mathrm{Op}\left(a\right) \right\|_{\xLtwo \to \xLtwo}
	\le C \left| a \right|_{S(1)}^{(k)} .
\end{equation}
\end{prpstn}

\begin{prpstn}[Sharp G{\aa}rding inequality --- {\cite[Proposition B.6]{P:23}}] \label{prop:Gaarding}
There exist a constant $c_d > 0$ and an integer $k_d \ge 0$ depending only on the dimension $d$ such that the following holds. For any real-valued symbol $a \in S(1)$ satisfying $a \ge 0$, one has
\begin{equation*}
\mathrm{Op}(a)
	\ge - c_d \left| \mathrm{Hess} \, a \right|_{S(1)}^{(k_d)} .
\end{equation*}
\end{prpstn}

%%-----------------------------
%%      your bibliography
%%-----------------------------

\small
\bibliographystyle{alpha}
\bibliography{biblio}

\end{document}